\documentclass[11pt]{amsart}
\usepackage{amscd, amsfonts, amsthm, graphicx, amssymb}
\title{Measure-preserving homeomorphisms of noncompact manifolds and mass flow toward ends} 
\author{Tatsuhiko Yagasaki}
\address{Division of Mathematics, Faculty of Engineering and Design, 
Kyoto Institute of Technology, 
Matsugasaki, Sakyoku, Kyoto 606-8585, Japan}
\email{yagasaki@kit.ac.jp}

\newtheorem{theorem}{Theorem}[section]
\newtheorem*{thm'}{Theorem 1.2$'$}
\newtheorem{proposition}{Proposition}[section] 
\newtheorem{corollary}{Corollary}[section] 
\newtheorem{lemma}{Lemma}[section]

\theoremstyle{definition}
\newtheorem{defi}{Definition}[section]

\def \cal {\mathcal}
\def \phi {\varphi}
\def \ds {\displaystyle}

\def \lra {\longrightarrow}

\def \e {\varepsilon}
\def \bs {\boldsymbol}

\begin{document}
\baselineskip 6 mm

\thispagestyle{empty}

\maketitle

\begin{abstract}
Suppose $M$ is a noncompact connected $n$-manifold and $\omega$ 
is a good Radon measure of $M$ with $\omega(\partial M) = 0$. 
Let ${\cal H}(M; \omega)$ denote the group of $\omega$-preserving homeomorphisms of $M$ 
equipped with the compact-open topology and ${\cal H}_E(M; \omega)$ 
denote the subgroup consisting of all $h \in {\cal H}(M; \omega)$ which fix the ends of $M$. 
 S.~R.~Alpern and V.\,S.\,Prasad introduced 
the topological vector space ${\cal S}(M, \omega)$ of end charges of $M$ and 
the end charge homomorphism $c^\omega : {\cal H}_E(M; \omega) \to {\cal S}(M, \omega)$, 
which measures for each $h \in {\cal H}_E(M; \omega)$ the mass flow toward ends induced by $h$. 
 We show that the map $c^\omega$ has a continuous section. 
This induces the factorization ${\cal H}_E(M; \omega) \cong {\rm Ker}\,c^\omega \times {\cal S}(M, \omega)$ 
and implies that ${\rm Ker}\,c^\omega$ is a strong deformation retract of ${\cal H}_E(M; \omega)$. 
\end{abstract}


\section{Introduction} 

This article is a continuation of the study of 
groups of measure-preserving homeomorphisms of noncompact topological manifolds \cite{Be1,Be2, Be3,Ya1}. 
Suppose $M$ is a noncompact connected $n$-manifold and $\omega$ 
is a good Radon measure of $M$ with $\omega(\partial M) = 0$. 
Let ${\cal H}(M; \omega)$ denote the group 
of $\omega$-preserving homeomorphisms of $M$ 
equipped with the compact-open topology.
In the study of this group, the space $E_M$ of ends of $M$ plays a significant role. 
Let $E_M^\omega$ denote the open subset of $E_M$ consisting of $\omega$-finite ends of $M$ and 
let ${\cal H}_{E_M}(M; \omega)$ denote the subgroup 
consisting of all $h \in {\cal H}(M; \omega)$ which fix the ends of $M$. 

In \cite{AP} S.~R.~Alpern and V.\,S.\,Prasad introduced the end charge homomorphism 
\[ c^\omega : {\cal H}_{E_M}(M; \omega) \lra {\cal S}(M, \omega). \] 
An end charge of $M$ is a finitely additive signed measure 
on the algebra of clopen subsets of $E_M$. 
Let ${\cal S}(E_M)$ denote the topological linear space of all end charges of $M$ with the weak topology and 
let ${\cal S}(M, \omega)$ denote the linear subspace of ${\cal S}(E_M)$ consisting of end charges $c$ of $M$ with $c(E_M) = 0$ and $c |_{E_M^\omega} = 0$. 
For each $h \in {\cal H}_{E_M}(M; \omega)$ an end charge $c^\omega_h \in {\cal S}(M, \omega)$ is defined by 
\[ c^\omega_h(E_C) = \omega(C - h(C)) - \omega(h(C) - C), \]
where $C$ is any Borel subset of $M$ such that ${\rm Fr}\,C$ is compact 
and $E_C \subset E_M$ is the set of ends of $C$. 
This quantity is the total $\omega$\,-\,volume (or mass) transfered by $h$ into $C$ and into $E_C$ in the last. 
Hence, the end charge $c^\omega_h$ measures mass flow toward ends induced by $h$. 

In \cite{Be3} R.~Berlanga showed that 
the group ${\cal H}(M; \omega)$ is a strong deformation retract of 
the group ${\cal H}(M; \omega\mbox{-e-reg})$ consisting of  
$\omega$-end-regular homeomorphisms of $M$. 
The group ${\cal H}(M; \omega\mbox{-e-reg})$ acts continuously on the space ${\cal M}_g^\partial(M, \omega\mbox{-e-reg})_{ew}^\ast$ of good Radon measures $\mu$ on $M$ such that $\mu(M) = \omega(M)$,  
$E_M^\mu = E_M^\omega$ and $\mu$ and $\omega$ have the same null sets, equipped with the finite-end weak topology. 
He showed that the orbit map $\pi : {\cal H}(M; \omega\mbox{-e-reg})
\lra {\cal M}_g^\partial(M, \omega\mbox{-e-reg})_{ew}$ : $h \longmapsto h_\ast \omega$ has a continuous section. 
This section induces the factorization 
${\cal H}(M; \omega\mbox{-e-reg}) \cong {\cal H}(M; \omega) \times {\cal M}_g^\partial(M, \omega\mbox{-e-reg})_{ew}^\ast$ and 
this yields the strong deformation retraction of ${\cal H}(M; \omega\mbox{-e-reg})$ onto ${\cal H}(M; \omega)$. 

In this article we use a similar strategy and investigate the internal structure of the group ${\cal H}(M; \omega)$. 
The group ${\cal H}_{E_M}(M, \omega)$ acts continuously on ${\cal S}(M, \omega)$ by 
$h \cdot a = c^\omega_h + a$ ($h \in {\cal H}_{E_M}(M, \omega)$, $a \in {\cal S}(M, \omega)$) 
and the end charge homomorphism $c^\omega : {\cal H}_{E_M}(M, \omega) \to {\cal S}(M, \omega)$ 
coincides with the orbit map at $0 \in {\cal S}(M, \omega)$. 
We extend the argument in \cite{Be3} and  show that 
the map $c^\omega$ admits a continuous (non-homomorphic) section. 

Suppose $M^n$ is a noncompact connected separable metrizable 
$n$-manifold and $\omega \in {\cal M}_g^\partial(M)$.  

\begin{theorem}
There exists a continuous map 
$s : {\cal S}(M, \omega) \to {\cal H}_\partial(M, \omega)_1$ 
such that $c^\omega s = id$ and $s(0) = id_M$. 
\end{theorem}

\begin{theorem} 
Suppose $P$ is any topological space and $\mu : P \to {\cal M}_g^\partial(M, \omega\mbox{-reg})$ and $a : P \to {\cal S}(E_M)$ 
are continuous maps such that $a_p \in {\cal S}(M; \mu_p)$ $(p \in P)$.  
Then there exists a continuous map 
$h : P \to {\cal H}_\partial(M, \omega\mbox{-reg})_1$ 
such that for each $p \in P$ \\
\hspace{15mm} 
$(1)$ $h_p \in {\cal H}_\partial(M, \mu_p)_1$, \hspace{4mm} 
$(2)$ $c^{\mu_p}_{h_p} = a_p$, \hspace{4mm}  
$(3)$ if $a_p = 0$, then $h_p = id_M$. 
\end{theorem}

Theorem 1.2 is a slight generalization of Theorem 1.1. 
The existence of a section for the map $c^\omega$ and 
the contractibility of the base space ${\cal S}(M, \omega)$ 
imply the following consequences.  

\begin{corollary} 
$(1)$ ${\cal H}_{E_M}(M; \omega) \cong {\rm Ker}\,c^\omega \times {\cal S}(M, \omega)$. 

$(2)$ ${\rm Ker}\,c^\omega$ is a strong deformation retract of ${\cal H}_{E_M}(M; \omega)$. 
\end{corollary}

The group ${\rm Ker}\,c^\omega$ contains the subgroup
${\cal H}^c(M; \omega)$ consisting of 
$\omega$-preserving homeomorphisms with compact support. 
The condition $c^\omega_h = 0$ means that 
any compact part of $h$ can be separated from the ``remaining part'' of $h$. 
From the argument in \cite{AP} it follows that 
for any $f \in {\rm Ker}\,c^\omega \cap {\cal H}(M)_1$ and any 
compact subset $K$ of $M$ there exists a compact connected $n$-submanifold $N$ of $M$ with $K \subset N$ and $h \in {\cal H}_{M - N}(M, \omega)_1$ with $h|_K = f|_K$. 
This implies that 
the subgroup ${\cal H}^c(M, \omega)_1^\ast$ is dense in ${\rm Ker}\,c^\omega \cap {\cal H}(M)_1$. 
In a succeeding work we will show that in $n =2$ 
the subgroup ${\cal H}^c(M, \omega)_1^\ast$ is homotopy dense in ${\rm Ker}\,c^\omega \cap {\cal H}(M)_1$.  
In \cite{Ya2} we have obtained some versions of Theorem 1.1 and \cite[Theorem 4.1]{Be3} for smooth manifolds and volume-preserving diffeomorphisms. 

This paper is organized as follows.  
Section 2 contains fundamentals on 
end compactifications, spaces of Radon measures 
and groups of measure-preserving homeomorphisms. 
Section 3 is devoted to basics on end charge homomorphisms and 
related notions. 
This section also includes generalities on morphisms induced from proper maps. 
Section 4 contains the proof of Theorem 1.2 in the cube case. 
The general case is treated in Section 5.


\section{Radon measures and end charge homomorphism} 

Throughout this section $X$ is a connected, locally connected, locally compact, separable metrizable space. 
We use the following notations : 
${\cal F}(X)$, ${\cal K}(X)$ and ${\cal C}(X)$  denote the sets of closed subsets, compact subsets, and connected components of $X$. 
${\cal B}(X)$ and ${\cal Q}(X)$ denote the $\sigma$-algebra of Borel subsets and 
the algebra of clopen subsets of $X$ respectively. 

When $A$ is a subset of $X$, the symbols 
${\rm Fr}_X A$, ${\rm cl}_X A$ and ${\rm Int}_X A$ 
denote the frontier, closure and interior of $A$ relative to $X$. 
When $M$ is a manifold, $\partial = \partial M$ and ${\rm Int}\,M$ 
denote the boundary and interior of $M$ as a manifold.

\subsection{Groups of homeomorphisms} 
For a space $X$ and a subset $A \subset X$ 
the symbol ${\cal H}_A(X)$ denotes 
the group of homeomorphisms $h$ of $X$ onto itself with $h|_A = id_A$ equipped with the compact-open topology. 
The group ${\cal H}_A(X)$ is a topological group 
(since $X$ is locally compact and locally connected). 

The support of $h \in {\cal H}(X)$ is defined by ${\rm Supp}\,h = cl_X\,\{ x \in X \mid h(x) \neq x \}$.
We set ${\cal H}_A^c(X) = \{ h \in {\cal H}_A(X) \mid {\rm Supp}\,h : \mbox{ compact}\}$. 
For any subgroup ${\cal G}$ of ${\cal H}(X)$, 
the symbol ${\cal G}_1$ denotes the path-component of $id_M$ in ${\cal G}$. 
When ${\cal G} \subset {\cal H}^c(X)$, by ${\cal G}_1^\ast$ we denote the subgroup of ${\cal G}_1$ 
consisting of $h \in {\cal G}$ which admits an isotopy $h_t \in {\cal G}$ ($t \in [0, 1]$) 
such that $h_0 = id_X$, $h_1 = h$ and there exists $K \in {\cal K}(X)$ with ${\rm Supp}\,h_t \subset K$ ($t \in [0,1]$).  

\subsection{End compactifications}  (cf.\,\cite{AP, Be3})
Suppose $X$ is a noncompact, connected, locally connected, locally compact, separable metrizable space. 
An end of $X$ is a function $e$ which assigns an $e(K) \in {\cal C}(X - K)$ to each $K \in {\cal K}(X)$ such that 
$e(K_1) \supset e(K_2)$ if $K_1 \subset K_2$. 
The set of ends of $X$ is denoted by $E_X$. 
The end compactification of $X$ is the space $\overline{X} = X \cup E_X$ 
equipped with the topology defined by the following conditions: 
(i) $X$ is an open subspace of $\overline{X}$, 
(ii) the fundamental open neighborhoods of $e \in E_X$ are given by 
\[ N(e, K) = e(K) \,\cup \,\{ e' \in E_X \mid e'(K) = e(K)\} \hspace{2mm} (K \in {\cal K}(X)). \] 
Then $\overline{X}$ is a connected, locally connected, compact, metrizable space, 
$X$ is a dense open subset of $\overline{X}$ and  
$E_X$ is a compact 0-dimensional subset of $\overline{X}$. 

Let ${\cal B}_c(X) = \{ C \in {\cal B}(X) \mid {\rm Fr}_X\,C : \mbox{compact} \}$.  
For each $C \in {\cal B}_c(X)$ let 
\[ E_C = \{ e \in E_X \mid e(K) \subset C \mbox{ for some } K \in {\cal K}(X)\}
\hspace{5mm} \text{and} \hspace{5mm} \overline{C} = C \cup E_C \subset \overline{X}. \]  
Then $E_C \in {\cal Q}(E_X)$ and $\overline{C}$ is a neighborhood of $E_C$ in $\overline{X}$ with $\overline{C} \cap E_X = E_C$.   
For $C, D \in {\cal B}_c(X)$, $E_C = E_D$ iff 
$C \Delta D = (C - D) \cup (D - C)$ is relatively compact (i.e., has the compact closure) in $X$. 

For $h \in {\cal H}(X)$ and $e \in E_X$ we define $h(e) \in E_X$ by $h(e)(K) = h(e(h^{-1}(K)))$ $(K \in {\cal K}(X))$. 
Each $h \in {\cal H}(X)$ has a unique extension $\overline{h} \in {\cal H}(\overline{X})$ defined by 
$\overline{h}(e) = h(e)$ $(e \in E_X)$. 
The map \ ${\cal H}(X) \to {\cal H}(\overline{X})$ : $h \mapsto \overline{h}$ \ is a continuous group homomorphism.  
We set ${\cal H}_{A \cup E_X}(X) = \{ h \in {\cal H}_A(X) \mid \overline{h}|_{E_X} = id_{E_X}\}$. 
Note that ${\cal H}_{A \cup E_X}(X)_1 = {\cal H}_{A}(X)_1$ and that 
if $C \in {\cal B}_c(X)$ and $h \in {\cal H}_{E_X}(X)$, then $h(C) \in {\cal B}_c(X)$ and $E_{h(C)} = E_C$. 


\subsection{Space of Radon measures} 

Next we recall general facts on spaces of Radon measures cf.\,\cite{AP, Be3, Fa}. 
Suppose $X$ is a connected, locally connected, locally compact, separable metrizable space. 
A {\em Radon measure} on $X$ is a measure $\mu$ on $(X, {\cal B}(X))$ 
such that $\mu(K) < \infty$ for any $K \in {\cal K}(X)$. 
A Radon measure $\mu$ on $X$ is said to be {\em good} if 
$\mu(p) = 0$ for any point $p \in X$ and $\mu(U) > 0$ for any nonempty open subset $U$ of $X$.  

Let ${\cal M}(X)$ denote the space of Radon measures $\mu$ on $X$  
equipped with  the {\em weak topology}. 
This topology is the weakest topology such that the function 
\[ \mbox{$\Phi_f : {\cal M}(X) \to {\Bbb R}$ : \ $\ds \Phi_f(\mu) = \int_X f \,d\mu$} \]
is continuous for any continuous function $f : X \to {\Bbb R}$ with compact support. 
Let ${\cal M}_g(X)$ denote the subspace of good Radon measures $\mu$ on $X$ and 
for $A \in {\cal B}(X)$ we set ${\cal M}^A(X) = \{ \mu \in {\cal M}(X) \mid \mu(A) = 0 \}$ and 
${\cal M}_g^A(X) = {\cal M}_g(X) \cap {\cal M}^A(X)$. 

For $\mu \in {\cal M}_g(X)$ and $A \in {\cal B}(X)$ 
the {\em restriction} $\mu|_A \in {\cal M}_g(A)$ 
is defined by $(\mu|_A)(B) = \mu(B)$ \ ($B \in {\cal B}(A)$). 

\begin{lemma} {\rm (\cite[Lemma 2.2]{Be3})} Let $A \in {\cal F}(X)$ and $K \in {\cal K}(X)$. 
\begin{itemize}
\item[(i)\,] The restriction map \  
${\cal M}^{{\rm Fr}A}(X) \lra {\cal M}(A)$ : $\mu \longmapsto \mu|_A$ \ is continuous.  
\item[(ii)] The evaluation map \ ${\cal M}^{{\rm Fr}K}(X) \lra {\Bbb R}$ : $\mu \longmapsto \mu(K)$ \ is continuous. 
\end{itemize} 
\end{lemma}

Let $\omega \in {\cal M}_g(X)$. 
We say that an end $e \in E_X$ is $\omega$-finite if $\omega(e(K)) < \infty$ for some $K \in {\cal K}(X)$. 
Let $E_X^\omega = \{ e \in E_X \mid e : \text{$\omega$-finite}\,\}$. This is an open subset of $E_X$ and 
for $C \in {\cal B}_c(X)$ we have 
$E_C \subset E_X^\omega$ iff $\omega(C) < \infty$. 

\begin{defi} (1) $\mu \in {\cal M}_g(M)$ is said to be 
\begin{itemize}
\item[(i)\,] {\em $\omega$-regular} if $\mu$ has the same null sets as $\omega$ (i.e., $\mu(B) = 0$ iff $\omega(B) = 0$ for any $B \in {\cal B}(X)$). 
\item[(ii)] {\em $\omega$-end-regular} if $\mu$ is $\omega$-regular and $E_M^\mu = E_M^\omega$. 
\end{itemize} 
(2) ${\cal M}^A_g(X, \omega\mbox{(-e)-reg}) = 
\big\{  \mu \in {\cal M}^A_g(X) \mid \mu : \omega\mbox{(-end)-regular} \big\}$ (the weak topology)
\end{defi}

The group ${\cal H}(X)$ acts continuously on ${\cal M}(X)$ 
by $h \cdot \mu = h_\ast \mu$, where $h_\ast \mu$ is defined by 
$(h_\ast \mu)(B) = \mu(h^{-1}(B))$ $(B \in {\cal B}(X))$. 

\begin{defi} (1) $h \in {\cal H}(X)$ is said to be 
\begin{itemize}
\item[(i)\ ] {\em $\omega$-preserving} if $h_\ast \omega = \omega$ \ \  
(i.e., $\omega(h(B)) = \omega(B)$ for any $B \in {\cal B}(X)$), 
\item[(ii)\,] {\em $\omega$-regular} if $h$ preserves $\omega$-null sets \ \ 
(i.e., $\omega(h(B)) = 0$ iff $\omega(B) = 0$ for any $B \in {\cal B}(X)$),  
\item[(iii)] {\em $\omega$-end-regular} if $h$ is $\omega$-regular and $\overline{h}(E_X^\omega) = E_X^\omega$. 
\end{itemize} 
(2) ${\cal H}(X; \omega) = \{ h \in {\cal H}(X) \mid h : \omega\mbox{-preserving}\}$, \  
${\cal H}(X; \omega\mbox{(-e)-reg}) = \{ h \in {\cal H}(X) \mid h : \omega\mbox{(-end)-regular}\}$
\end{defi}

Suppose $M$ is a compact connected $n$-manifold.
The von Neumann-Oxtoby-Ulam theorem \cite{OU} asserts that 
if $\mu, \nu \in {\cal M}_g^\partial(M)$ and $\mu(M) = \nu(M)$, then 
there exists $h \in {\cal H}_\partial(M)_1$ such that $h_\ast \mu = \nu$.
A.\,Fathi \cite{Fa} obtained a parameter version of this theorem. 

\begin{theorem} Suppose $M$ is a compact connected $n$-manifold 
and $\omega \in {\cal M}_g^\partial(M)$. 
Suppose $\mu, \nu : P \to {\cal M}_g^\partial(M; \omega\mbox{-reg})$ 
are continuous maps with 
$\mu_p(M) = \nu_p(M)$ $(p \in P)$. 
Then there exists a continuous map 
$h : P \to {\cal H}_\partial(M; \omega\mbox{-reg})_1$ such that for each $p \in P$ 
{\rm (i)} $(h_p)_\ast \mu_p = \nu_p$ and {\rm (ii)} if $\mu_p = \nu_p$ then $h_p = id_M$.  
\end{theorem}

In \cite{Be3} R.\,Berlanga  obtained 
a similar theorem for a noncompact connected $n$-manifold $M$. 
We use the following consequence of \cite[Proposition 5.1\,(2)]{Be3}. 

\begin{lemma} Suppose $M$ is a noncompact connected $n$-manifold and $\omega \in {\cal M}_g^\partial(M)$. Then 
we have 
${\cal H}_\partial(M; \omega) \cap {\cal H}_\partial(M; \omega\mbox{-reg})_1 
= {\cal H}_\partial(M; \omega)_1$. 
\end{lemma}


\section{End charge homomorphism} 

\subsection{End charge homomorphism} 
We recall basic properties of the end charge homomorphism defined in \cite[Section 14]{AP}. 
Suppose $X$ is a connected, locally connected, locally compact, separable, metrizable space and 
$\omega \in {\cal M}(X)$. 

An {\em end charge} of $X$ is a finitely additive signed measure $c$ on ${\cal Q}(E_X)$, that is, 
a function $c : {\cal Q}(E_X) \to {\Bbb R}$ which satisfies the following condition:  
\[ \mbox{$c(F \cup G) = c(F) + c(G)$ \ for \ $F, G \in {\cal Q}(E_X)$ \ with \ $F \cap G = \emptyset$.} \]  
Let ${\cal S}(E_X)$ denote the space of end 
charges $c$ of $X$ with the {\em weak topology} (or the product topology). 
This topology is the weakest topology such that the function 
\[ \Psi_F : {\cal S}(E_X) \lra {\Bbb R}  : \ \Psi_F(c) = c(F) \] 
is continuous for any $F \in {\cal Q}(E_X)$. 
For a subset $U \subset E_X$ let 
\[ \mbox{${\cal S}_0(E_X, U) = \big\{ c \in {\cal S}(E_X) \mid$ (i) $c(F) = 0$ for 
$F \in {\cal Q}(E_X)$ with $F \subset U$ and (ii) $c(E_X) = 0 \,\big\}$} \] 
(with the weak topology). 
Then ${\cal S}(E_X)$ is a topological linear space and ${\cal S}_0(E_X, U)$ is a linear subspace. 
For $\omega \in {\cal M}(X)$ we set ${\cal S}(X, \omega) = {\cal S}_0(E_X, E_X^\omega)$. 

For $h \in {\cal H}_{E_X}(X, \omega)$ 
the end charge $c_h^\omega \in {\cal S}(X, \omega)$ is defined as follows: 
For any $F \in {\cal Q}(E_X)$ there exists $C \in {\cal B}_c(X)$ with $E_C = F$.  
Since $\overline{h}|_{E_X} = id$, it follows that $E_C = E_{h(C)}$ and that 
$C \Delta \,h(C)$ is relatively compact in $X$. 
Thus $\omega(C - h(C)), \omega(h(C) - C) < \infty$ and we can define as  
\[ c_h^\omega(F) = \omega(C - h(C)) - \omega(h(C) - C) \in {\Bbb R}.  \] 
This quantity is independent of the choice of $C$. 

\begin{proposition} The map $c^\omega : {\cal H}_{E_X}(X, \omega) \lra {\cal S}(X, \omega)$ is   
a continuous group homomorphism $($\cite[Section 14.9, Lemma 14.21\,(iv)]{AP}$)$. 
\end{proposition} 

\subsection{Related notions} 
In the proof of Theorem 1.2 it is necessary to measure volumes transfered into various regions by homeomorphisms (which are not measure-preserving). 
For this purpose we introduce some notations. 

For $A, B \in {\cal B}(X)$ we write $A \sim_c B$ if $A \Delta B$ is relatively compact in $X$. This is an equivalence relation and 
for $A, B \in {\cal B}_c(X)$ we have (i) $A \sim_c B$ iff $E_A = E_B$ and 
(ii) $A \sim_c h(A)$ for any $h \in {\cal H}_{E_X}(X)$. 

Similarly, for $\mu \in {\cal M}(X)$ and $A, B \in {\cal B}(X)$ we write $A \sim_\mu B$ if $\mu(A \Delta B) < \infty$. 
This is also an equivalence relation and $A  \sim_c B$ implies $A \sim_\mu B$. If $A \sim_\mu B$, then  
we can consider the following quantity: 
\[ J^\mu(A, B) = \mu(A - B) - \mu(B - A) \ \in \ {\Bbb R}. \]

This measures the difference of $\mu$-volumes of $A$ and $B$ when 
$A$ and $B$ differ only in a finite volume part. 
If $C \in {\cal B}_c(X)$ and $h \in {\cal H}_{E_X}(X)$, then 
$J^\mu(h^{-1}(C), C)$ is just the total $\mu$\,-\,mass transfered  into $C$ by $h$. If $h \in {\cal H}_{E_X}(X, \mu)$, then 
$J^\mu(h^{-1}(C), C) = J^\mu(C, h(C)) = c_h^\mu(E_C)$. 

This quantity has the following formal properties: 

\begin{lemma} \label{l-difference of volume} 
Suppose $\mu \in {\cal M}(X)$ and $A, B, C, D \in {\cal B}(X)$.

$(1)$ If $A \sim_\mu B$ and $\mu(A) < \infty$, then $\mu(B) < \infty$ and $J^\mu(A, B) = \mu(A) - \mu(B)$. 

$(2)$ If $A \sim_\mu B \sim_\mu C$, then 
$J^\mu(A, B) + J^\mu(B, C) = J^\mu(A, C)$. 

$(3)$ If $A \sim_\mu C$, $B \sim_\mu D$, then 
\begin{itemize}
\item[(i)\,] $A \cup B \sim_\mu C \cup D$ since $(A \cup B) \Delta (C \cup D) \subset (A \Delta C) \cup (B \Delta D)$, 
\item[(ii)] if $A \cap B  = C \cap D = \emptyset$, then $J^\mu(A \cup B, C \cup D) = J^\mu(A, C) + J^\mu(B, D)$.
\end{itemize} 

$(4)$ If $h \in {\cal H}(X)$ and $A \sim_{h_\ast \mu} B$, then 
$h^{-1}(A) \sim_\mu h^{-1}(B)$ and 
$J^{h_\ast \mu}(A, B) = J^\mu(h^{-1}(A), h^{-1}(B))$. 
\end{lemma} 

\begin{lemma} Suppose $\omega \in {\cal M}(X)$ and $A, B \in {\cal B}_c(X)$, $A \sim_c B$, $\omega({\rm Fr}\,A) = \omega({\rm Fr}\,B) = 0$. Then the function 
$$\Phi : {\cal M}(X: \omega\mbox{-reg}) \times {\cal H}_{E_X}(X; \omega\mbox{-reg})^2 \lra {\Bbb R} : \ \ \Phi(\mu, f, g) =  J^\mu(f(A), g(B))$$ is continuous. 
\end{lemma} 

\begin{proof} Since 
$J^\mu(f(A), g(B)) = J^\mu(f(A), A) + J^\mu(A, B) + J^\mu(B, g(B))$ and \\
$\mu(A - f(A)) = (f^{-1}_\ast \mu)(f^{-1}(A) - A)$, it suffices to verify the continuity of the following function:  
$${\cal M}(X: \omega\mbox{-reg}) \times {\cal H}_{E_X}(X; \omega\mbox{-reg}) \lra {\Bbb R} : \ \ (\mu, f) \longmapsto  \mu(f(A) - A).$$ 

Given $(\mu, f)$ and $\e > 0$. 
Since $f({\rm Fr}\,A)$ is a compact $\mu$-null set, it has a compact neighborhood $K$ such that $\mu(K) < \e$.
There exists a neighborhood ${\cal U}$ of $f$ in ${\cal H}_{E_X}(X, \omega\mbox{-reg})$ such that $f(A) \Delta g(A) \subset K$ ($g \in {\cal U}$). 

The function $\nu(f(A) - A)$ is continous in $\nu$. 
In fact,  
${\rm Fr}\,(f(A) - A) \subset {\rm Fr}\,A \cup {\rm Fr}\,f(A)$ and the latter is a $\nu$-null set since $\nu$ is $\omega$-regular. 
Thus, we have $\nu({\rm Fr}\,(f(A) - A)) = 0$ and the claim follows from Lemma 2.1\,(ii). 
Also note that the function ${\cal M}(X) \to {\Bbb R} : \nu \mapsto \nu(K)$ is upper semi-continuous (\cite[Lemma 2.1]{Be3}). 
Therefore, there exists a neighborhood ${\cal V}$ of $\mu$ in ${\cal M}(X ; \omega\mbox{-reg})$ such that
$$|\nu(f(A) - A) - \mu(f(A) - A)| < \e \ \ \text{and} \ \ \nu(K) < \e \ \  (\nu \in {\cal V}).$$ 
Take any $(\nu, g) \in {\cal V} \times {\cal U}$. 
Since 
$(f(A) - A) \Delta (g(A) - A) \subset f(A) \Delta g(A) \subset K$, we have 
$$|\nu(g(A) - A) - \nu(f(A) - A)| \leq \nu(K) < \e.$$
(In general, $|\nu(A) - \nu(B)| \leq \nu(A \Delta B)$.) It follows that 
$|\nu(g(A) - A) - \mu(f(A) - A)| < 2\e.$
\end{proof} 

According to \cite{Be3} we say that 
continuous maps $\mu, \nu : P \to {\cal M}(X)$ 
are {\em compactly related} and write $\mu \sim_c \nu$ 
if each $p \in P$ admits a neighborhood $U$ in $P$ and $K_p \in {\cal K}(X)$ such that 
$\mu_q = \nu_q$ on $M - K_p$ $(q \in U)$. 
(If $P$ is a singleton, this is just a condition on $\mu, \nu \in {\cal M}(X)$.) 
This is an equivalence relation and  
if $\mu \sim_c \nu$, then for any $C \in {\cal B}(X)$ we can define a function 
$(\mu - \nu)(C) : P \to {\Bbb R}$ by 
\[ (\mu - \nu)(C)_p = \mu_p(C \cap K_p) - \nu_p(C \cap K_p). \] 
This definition is independent of the choice of $K_p$.
If $\omega \in {\cal M}(X)$, $\mu, \nu : P \to {\cal M}(X, \omega\mbox{-reg})$, 
$\mu \sim_c \nu$ 
and $C \in {\cal B}(X)$, $\omega({\rm Fr}\,C) = 0$, then 
the function $(\mu - \nu)(C) : P \to {\Bbb R}$  
is continuous. 

Suppose a continuous map 
$h : P \to {\cal H}^c(X)$ has {\em locally common compact support} 
(i.e., for each $p \in P$ there exists a neighborhood $U$ of $p$ in $P$ and $K \in {\cal K}(X)$ such that 
${\rm Supp}\,h_q \subset K$ $(q \in U)$). 
Then, $\mu \sim_c h_\ast \mu$ 
for any continuous map $\mu : P \to {\cal M}(X)$.

If $\mu \in {\cal M}(X)$, $A \in {\cal B}(X)$ and 
$f, g \in {\cal H}^c(X)$, then we have the following relation
\[ (f_\ast \mu - g_\ast \mu)(A) = J^\mu(f^{-1}(A), g^{-1}(A)). \] 

In the proof of Theorem 1.2 we use the quantity of the form $J^\mu(f^{-1}(A), g^{-1}(A))$
frequently. 
The above consideration 
means that this quantity can be translated to a  
quantity prescribed in term of measures and that 
the calculations on this quantity in the proof of Theorem 1.2 and 
the statements in Lemma \ref{l-difference of volume} reduce to 
the calculations and some ordinary properties on measures.  
However, the quantity $J^\mu(f^{-1}(A), g^{-1}(A))$ has an advantage that 
it is defined for $A \in {\cal B}_c(X)$ and $f, g \in {\cal H}_{E_X}(X)$. 
For example, we can take $f$ and $g$ as the limits of sequences $f_k, g_k \in {\cal H}^c(X)$. This fits our situation. 

\subsection{Morphisms induced from proper maps} 

Suppose $X$ and $Y$ are connected, locally connected, locally compact separable metrizable spaces and 
$f: X \to Y$ is a {\em proper} continuous map ($f^{-1}(K)$  is compact for any $K \in {\cal K}(Y)$). 
The map $f$ induces various continuous morphisms. 

(1) $f_\ast : {\cal M}(X) \to {\cal M}(Y)$ : 
For $\mu \in {\cal M}(X)$ the induced measure $f_\ast \mu \in {\cal M}(Y)$ is defined by $(f_\ast \mu)(B) = \mu(f^{-1}(B))$ ($B \in {\cal B}(Y)$). The map $f_\ast$ is continuous. 
If $A, B \in {\cal B}(Y)$ and $A \sim_{f_\ast \mu} B$, then 
$J^{f_\ast \mu}(A, B) = J^\mu(f^{-1}(A), f^{-1}(B))$ 
(cf. Lemma \ref{l-difference of volume}\,(4)).
\vskip 1mm 
(2) $\overline{f} : \overline{X} \to \overline{Y}$ : This is the unique continuous extension of $f$.  
For each $e \in E_X$ the end $f(e) \in E_Y$ is defined by assigning to each $K \in {\cal K}(Y)$ 
the unique component $f(e)(K) \in {\cal C}(Y - K)$ which contains $f(e(f^{-1}(K)))$.
The map $\overline{f}$ is defined by $\overline{f}(e) = f(e)$ ($e \in E_X$). 
For any $C \in {\cal B}_c(Y)$ we have $f^{-1}(C) \in {\cal B}_c(X)$ and $E_{f^{-1}(C)} = \overline{f}^{-1}(E_C)$. 
\vskip 1mm 
(3) $\overline{f}_\ast : {\cal S}(E_X) \to  {\cal S}(E_Y)$ : 
This is a continuous linear map induced from the map $\overline{f} : E_X \to E_Y$. 
\break For each $c \in {\cal S}(E_X)$ the end charge 
$\overline{f}_\ast c \in {\cal S}(E_Y)$ is defined by 
$(\overline{f}_\ast c)(F) = c(\overline{f}^{-1}(F))$ \break $(F \in {\cal Q}(E_Y))$.
It induces the restriction 
$\overline{f}_\ast : {\cal S}_0(E_X, U) \to {\cal S}_0(E_Y, V)$  
for any $V \subset E_Y$ and $U \subset E_X$ with $\overline{f}^{-1}(V) \subset U$.
Let $\omega \in {\cal M}(X)$. 
Since $\overline{f}^{-1}(E_Y^{f_\ast \omega}) \subset E_X^\omega$, 
we obtain the restriction $\overline{f}_\ast : {\cal S}(X, \omega) \to {\cal S}(Y, f_\ast \omega)$. 
If $\overline{f} : E_X \to E_Y$ is injective, then 
$\overline{f}^{-1}(E_Y^{f_\ast \omega}) = E_X^\omega$. Therefore, if $\overline{f} : E_X \to E_Y$ is a homeomorphism, then 
$\overline{f}_\ast : {\cal S}(X, \omega) \to {\cal S}(Y, f_\ast \omega)$ is also a  homeomorphism. 
\vskip 2mm 
Below we assume that the map $f : X \to Y$ satisfies the following additional conditions: \\
\hspace{10mm} 
\begin{tabular}[t]{cl}
$(\ast)_1$ & $C \in {\cal F}(X)$, ${\rm Int}_X\,C = \emptyset$ and $D \in {\cal F}(Y)$, \\[1.5mm] 
$(\ast)_2$ & $f(C) = D$ and $f$ maps $X - C$ homeomorphically onto $Y - D$. 
\end{tabular}
\vskip 3mm 
(4) $f^\ast : {\cal M}^D(Y) \to {\cal M}^C(X)$ : 
For each $\nu \in {\cal M}^D(Y)$ the measure $f^\ast \nu \in {\cal M}^C(X)$ is defined by 
$(f^\ast \nu)(B) = \nu(f(B - C)) \ \ (B \in {\cal B}(X))$. 
The map $f^\ast$ is a homeomorphism, whose inverse is  
the map $f_\ast : {\cal M}^C(X) \to {\cal M}^D(Y)$.  
For any $\omega \in {\cal M}_g^D(Y)$ these maps induce the reciprocal homeomorphisms \\[1mm] 
\hspace{12mm} $f_\ast : {\cal M}_g^C(X; f^\ast \omega\mbox{-reg}) \to {\cal M}_g^D(Y; \omega\mbox{-reg})$, \hspace{5mm} 
$f^\ast : {\cal M}_g^D(Y; \omega\mbox{-reg}) \to {\cal M}_g^C(X; f^\ast \omega\mbox{-reg})$. 
\vskip 1mm 
(5) $f_\ast : {\cal H}_C(X) \to {\cal H}_D(Y)$ : 
For each $h \in {\cal H}_C(X)$ there exists a unique $\underline{h} \in 
{\cal H}_D(Y)$ such that $\underline{h} f = f h$. 
The map $f_\ast$ is defined by $f_\ast h = \underline{h}$. 
This map is a continuous injection and induces the restrictions 
$f_\ast : {\cal H}_{C \cup E_X}(X) \to {\cal H}_{D \cup E_Y}(Y)$ 
and 
$f_\ast : {\cal H}_C(X, f^\ast \omega\mbox{-reg}) \to {\cal H}_D(Y, \omega\mbox{-reg})$,  
$f_\ast : {\cal H}_C(X, f^\ast \omega) \to {\cal H}_D(Y, \omega)$ 
for any $\omega \in {\cal M}^D(Y)$. 

\begin{lemma}\label{l-morphism} 
 Under the condition $(\ast)$, for any $\omega \in {\cal M}^D(Y)$ we have 
the following commutative diagram : \\[-2mm] 
\hspace{50mm} 
$\begin{array}[t]{cccc}
& c^{f^\ast \omega} & & \\[-2mm] 
{\cal H}_{C \cup E_X}(X, f^\ast \omega) & \lra & {\cal S}(X, f^\ast \omega) & \\[2mm] 
\Big\downarrow & \hspace{-48mm} f_\ast & \Big\downarrow  & \hspace{-15.5mm} \overline{f}_\ast \\[3mm]  
{\cal H}_{D \cup E_Y}(Y, \omega) & \lra & {\cal S}(Y, \omega).  & \\[-1mm] 
& c^\omega & & 
\end{array}$ 
\end{lemma} 


\section{Proof of Theorem 1.2 in the cube case} 

In this section we prove Theorem 1.2 in the cube case. 
According to \cite[Section 4]{Be3} we use the following notations: 
$I = [0,1]$, $I^n$ is the $n$-fold product of $I$,  
$I_1 = [1/3, 2/3] \times \{ (1/2,\cdots, 1/2, 1) \} \subset I^n$,  
$m$ is the Lebesgue measure on ${\Bbb R}^n$, $d$ is the standard Euclidean distance in ${\Bbb R}^n$ ($d(\bs{x}, \bs{y}) = \| \bs{x} - \bs{y}\|$), 
$E$ is a 0-dim compact subset of $\partial I^n$ ($E \subset I_1$ for $n \geq 2$), $M_0 = I^n - E$ and $m_0 = m|_{M_0}$. 
The pair $(\overline{M_0}, E_{M_0})$ is canonically identified with $(I^n, E)$. 
An {\em $n$-cubic balloon} in $I^n$ is a cube $A$ of the form $[0, \alpha]^n + \bs{v}$ for some $\alpha > 0$ and $\bs{v} \in {\Bbb R}^n$ such that $A \subset I^n$ and $A \cap \partial I^n = ([0,\alpha]^{n-1} \times \{ \alpha \}) + \bs{v}$. 
Let ${\cal D}(M_0)$ denote the set of PL $n$-disks $K$ in $M_0$ such that $cl_{I^n}(M_0 - K)$ is a finite disjoint union of $n$-cubic balloons $A$ in $I^n$ with $A \cap E \neq \emptyset$. For convenience, we add the emptyset $\emptyset$ as a member of ${\cal D}(M_0)$. 

\begin{thm'} Suppose 
$\mu : P \to {\cal M}_g^\partial(M_0, m_0\mbox{-reg})$ 
and $a : P \to {\cal S}(E_{M_0})$ 
are continuous maps such that $a_p \in {\cal S}(M_0, \mu_p)$ $(p \in P)$. 
Then there exists a continuous map 
$h : P \to {\cal H}_\partial(M_0, m_0\mbox{-reg})_1$ 
such that for each $p \in P$ \\
\hspace{15mm} 
$(1)$ $h_p \in {\cal H}_\partial(M_0, \mu_p)_1$, \hspace{6mm} 
$(2)$ $c^{\mu_p}_{h_p} = a_p$, 
\hspace{6mm} 
$(3)$ if $a_p = 0$, then $h_p = id_{M_0}$. 
\end{thm'} 

Theorem 1.2$'$ is proved in a series of lemmas. 
For the sake of notational simplicity, 
we write $f_\ast \mu = g_\ast \mu$ and  
$J^\mu(f(A), g(A)) = a(E_A)$ instead of 
${f_p}_\ast \mu_p = {g_p}_\ast \mu_p$ $(p \in P)$ and 
$J^{\mu_p}(f_p(A), g_p(A)) = a_p(E_A)$ $(p \in P)$. 

Below we assume that 
$\mu : P \to {\cal M}_g^\partial(M_0, m_0\mbox{-reg})$ 
and $a : P \to {\cal S}(E_{M_0})$ 
are continuous maps such that $a_p \in {\cal S}(M_0, \mu_p)$ $(p \in P)$. 
We consider the case $n \geq 2$. (The modification for $n = 1$ is obvious.) 

\begin{lemma}\label{l-basic}  Suppose $K, L \in {\cal D}(M_0)$, $K \subset {\rm Int}_{M_0}\,L$  and 
$f, g : P \to {\cal H}_{\partial}(M_0, m_0\mbox{-reg})_1$
are continuous maps such that \\[1mm] 
\hspace{10mm} {\rm (i)} $f_\ast \mu = g_\ast \mu$ on $K$, \hspace{4mm} {\rm (ii)} $J^\mu(f^{-1}(A), g^{-1}(A)) = a(E_A)$ \ $(A \in {\cal C}(cl_{M_0}(M_0 - K)))$. \\[1mm] 
Then there exists a continuous map $h : P \to {\cal H}^c_{\partial \cup K}(M_0, m_0\mbox{-reg})_1^\ast$ such that 
\begin{itemize} 
\item[(1)] $(hf)_\ast \mu = g_\ast \mu$ on $L$, 
\item[(2)] $J^\mu((hf)^{-1}(B), g^{-1}(B)) = a(E_B)$ $(B \in {\cal C}(cl_{M_0}(M_0 - L)))$, 
\item[(3)] $\big\{ h_p^{-1} \big\}_{p \in P}$ is equi-continuous on $cl_{M_0}(M_0 - L)$ with respect to $d|_{M_0}$, 
\item[(4)] if $p \in P$,  $a_p = 0$ and $f_p = g_p = id_{M_0}$, then $h_p = id_{M_0}$. 
\end{itemize} 
\end{lemma} 

\begin{proof} 
For each $A \in {\cal C}(cl_{M_0}(M - K))$ we construct  
a continuous map $\ell = \ell_A : P \to {\cal H}_\partial^c(A, m|_A\mbox{-reg})_1^\ast$ such that 
\begin{itemize}
\item[(1)$'$] $\ell_\ast ((f_\ast \mu)|_A) = g_\ast \mu$ on $A \cap L$,  
\item[(2)$'$] $J^\mu(f^{-1}\ell^{-1}(B), g^{-1}(B)) = a(E_B)$ 
$(B \in {\cal C}(cl_{M_0}(A - L)))$, 
\item[(3)$'$] $\big\{ \ell_p^{-1} \big\}_{p \in P}$ is equi-continuous on $cl_{M_0}(A - L)$ with respect to $d|_{A}$, 
\item[(4)$'$]  if $p \in P$, $a_p = 0$ and $f_p = g_p = id_{M_0}$, then $\ell_p = id_A$. 
\end{itemize} 
\noindent 
Then the map $h$ is defined by $h|_K = id_K$ and $h|_{A} = \ell_A$ 
$(A \in {\cal C}(cl_{M_0}(M_0 - K)))$.  

The map $\ell = \ell_A$ is constructed as follows. 
Let ${\cal C}(cl_{M_0}(A - L)) = \{ B_1, \cdots, B_m \}$. 
This is a disjoint family of $n$-cubic balloons with ends and 
we have $A = (A \cap L) \cup \big(\cup_{k=1}^m B_k \big)$ and $E_A = \cup_{k=1}^m E_{B_k}$. 
Set $N_k = (A \cap L) \cup \big(\cup_{i=k}^m B_i\big)$ ($k = 1, \cdots, m$) and $N_{m+1} = A \cap L$. 

We inductively construct continuous maps 
$\ell^k : P \to {\cal H}_{\partial}^c(A, m|_A\mbox{-reg})_1^\ast$  ($k = 1, \cdots, m$) such that  
\begin{itemize}
\item[$(2_k)$] $J^\mu\big(f^{-1}(\ell^k)^{-1}(B_j), g^{-1}(B_j)\big) = a(E_{B_j})$ 
$(j=1, \cdots, k)$,  
\item[$(3_k)$] $\big\{ (\ell^k_p)^{-1} \big\}_p$ is equi-continuous with respect to $d|_A$, 
\item[$(4_k)$] if $p \in P$, $a_p = 0$ and $f_p = g_p = id_{M_0}$, then $\ell^k_p = id_{A}$. 
\end{itemize} 

Suppose $\ell^{k-1}$ has been constructed. (For $k=1$ we put $\ell^0 \equiv id_{M_0}$.) 
Consider the PL $n$-disk $\overline{N_k} = \overline{B_k} \cup \overline{N_{k+1}}$ (recall that 
$\overline{N_k} = N_k \cup E_{N_k}  = cl_{I^n}N_k$). 
Since $\overline{B_k} \cap \overline{N_{k+1}}$ is a PL $(n-1)$-disk, 
we can find a one-parameter family of PL-maps $\phi_t : \overline{N_k} \to \overline{N_k}$ ($t \in [-1,1]$) such that 
\begin{itemize} 
\item[(a)] $\phi_0 = id$, $\phi_{1}(\overline{B_k}) = \overline{N_k}$,  
$\phi_{-1}(\overline{N_{k+1}}) = \overline{N_k}$ and $\phi_t = id$ on $\partial \overline{N_k}$ ($t \in [-1,1]$), 
\vskip 1mm 
\item[(b)] $\phi_t|_{N_k}$ ($t \in (-1,1)$) is an isotopy on $N_k$, $\phi_s(B_k) \subsetneqq \phi_t(B_k)$ $(-1 \leq s < t \leq 1)$ 
and \\ 
$\phi_t|_{N_k}$ ($t \in (-1,1)$) has locally common compact support.
\end{itemize} 
The map $\phi_t$ is obtained 
by enlarging $\overline{B_k}$ for $t \geq 0$ (engulfing $\overline{N_k}$ at $t = 1$) and 
shrinking $\overline{B_k}$ for $t \leq 0$ (collapsing at $t = -1$). 
The family $\phi_t$ ($t \in [-1, 1]$) is equi-continuous with respect to $d$, since it is a compact family. Thus $\phi_t|_{N_k}$ ($t \in (-1, 1)$) is also equi-continuous with respect to $d$. 
The maps $\phi_t$ ($t \in (-1,1)$) are $m$-regular since any PL-homeomorphism between subpolyhedra in ${\Bbb R}^n$ is $m$-regular. 

The map $\ell^k$ is defined as $\ell^k = \psi \,\ell^{k-1}$, where   
$\psi : P \to {\cal H}_{\partial \cup B_1 \cup \cdots \cup B_{k-1}}^c(A, m|_A\mbox{-reg})_1^\ast$ is defined by 
\[ \mbox{$\psi_p = \phi_{t(p)}^{\ \ \ -1}$ \ on \ $N_k$ \ \ and \ \ 
$\psi_p = id$ \ on \ $B_1 \cup \cdots \cup B_{k-1}$.} \] 
The parameter function $t = t(p) : P \to (-1, 1)$ is determined by the condition $(2_k)$ ($j = k$). We set $$\sigma^{k-1}_p \equiv {\ell^{k-1}_p}_\ast (({f_p}_\ast\mu_p)|_A) \in {\cal M}_g^\partial(A, m|_A\mbox{-reg}).$$ 
Then the identity for $j = k$ in the condition $(2_k)$ is equivalent to : \\[3mm] 
\hspace{3mm} $\begin{array}[c]{lcl} 
\multicolumn{3}{l}{a_p(E_{B_k}) - J^{\mu_p}\big(f_p^{-1}(\ell^{k-1}_p)^{-1}(B_k) , g_p^{-1}(B_k)\big)} \\[5mm]  
\hspace{8mm} 
&=& J^{\mu_p}\big(f_p^{-1}(\ell^k_p)^{-1}(B_k) , f_p^{-1}(\ell^{k-1}_p)^{-1}(B_k)\big)  
\ = \ J^{\mu_p}\big(f_p^{-1}(\ell^{k-1}_p)^{-1}\phi_t(B_k) , f_p^{-1}(\ell^{k-1}_p)^{-1}(B_k)\big) \\[4mm]
&=& J^{\sigma^{k-1}_p}\big(\phi_t(B_k), B_k\big) 
\ = \ 
\left\{ \hspace{-1mm} 
\begin{array}[c]{rll}
\sigma^{k-1}_p (\phi_t(B_k) - B_k) & \hspace{-2mm} \in \ [0, \sigma^{k-1}_p(N_{k+1})) & (t \in [0, 1)) \\[4mm]
- \sigma^{k-1}_p (B_k - \phi_t(B_k)) & \hspace{-2mm} \in \ (- \sigma^{k-1}_p(B_k), 0] & (t \in (-1, 0]).
\end{array}\right. 
 \end{array}$ \\[5mm] 
(Note that $a_p(E_{B_k}) = 0$ does not imply $t(p) = 0$.) 
This equation in $t$ is uniquely solved, 
once we check the next inequality: 
$$a_p(E_{B_k}) - J^{\mu_p}\big(f_p^{-1}(\ell^{k-1}_p)^{-1}(B_k) , g_p^{-1}(B_k)\big) 
\in \big(- \sigma^{k-1}_p(B_k), 
\sigma^{k-1}_p(N_{k+1})\big).$$ 
This is verified by the following observations: 
\vskip 1mm 
If $\sigma^{k-1}_p(B_k) = \mu_p\big(f_p^{-1}(\ell^{k-1}_p)^{-1}(B_k)\big) < \infty$, then 
$\mu_p(B_k) < \infty$ and $a_p(E_{B_k}) = 0$ since $a_p \in S(M_0, \mu_p)$. Thus \\[2mm] 
\hspace{10mm} 
$\begin{array}[t]{ccl}
\multicolumn{3}{l}{a_p(E_{B_k}) - J^{\mu_p}\big(f_p^{-1}(\ell^{k-1}_p)^{-1}(B_k) , g_p^{-1}(B_k)\big)} \\[2mm] 
\hspace{14mm} &=& - \Big( \mu_p\big(f_p^{-1}(\ell^{k-1}_p)^{-1}(B_k)\big) - 
\mu_p(g_p^{-1}(B_k))\Big) \ > \ -\sigma^{k-1}_p(B_k). 
\end{array}$ 
\vskip 4mm 
If $\sigma^{k-1}_p(N_{k+1}) = \mu_p\big(f_p^{-1}(\ell^{k-1}_p)^{-1}(N_{k+1})\big)< \infty$, then $\mu_p(N_{k+1}) < \infty$ and 
$a_p(E_{B_j}) = 0$ ($j=k+1,$ $\cdots, m$). Since \\[2mm] 
\hspace{15mm} 
$\begin{array}[t]{l}
\ds \sum_{j=1}^{m} a(E_{B_j}) = a(E_A) = J^{\mu}(f^{-1}(A), g^{-1}(A)), \hspace{10mm} 
A \ =\ (\ell^{k-1})^{-1}(A), \\[6mm] 
a(E_{B_j}) \ = \ J^{\mu}\big(f^{-1}(\ell^{k-1})^{-1}(B_j) , g^{-1}(B_j)\big) 
\hspace{5mm} (j=1, \cdots, k-1), 
\end{array}$ \\[3mm] 
it follows that  \\ 
\hspace{15mm} 
$\begin{array}[t]{rcl} 
a_p(E_{B_k}) 
&=& \ds \sum_{j=1}^m a_p(E_{B_j}) - \left(\sum_{j=1}^{k-1} a_p(E_{B_j}) + \sum_{j=k+1}^m a_p(E_{B_j}) \right) \\[6mm]  
& = & \ds J^{\mu_p}\big(f_p^{-1}(\ell^{k-1}_p)^{-1}(A), g_p^{-1}(A)\big) 
 - \sum_{j=1}^{k-1} J^{\mu_p}\big(f_p^{-1}(\ell^{k-1}_p)^{-1}(B_j), g_p^{-1}(B_j)\big) \\[7mm] 
& = & J^{\mu_p}\big(f_p^{-1}(\ell^{k-1}_p)^{-1}(N_k) , g_p^{-1}(N_k)\big). 
\end{array}$ \\[7mm]  
\hspace{15mm} 
$\begin{array}[t]{lcl} 
\multicolumn{3}{l}{a_p(E_{B_k}) - J^{\mu_p}\big(f_p^{-1}(\ell^{k-1}_p)^{-1}(B_k) , g_p^{-1}(B_k)\big) \ = \ J^{\mu_p}\big(f_p^{-1}(\ell^{k-1}_p)^{-1}(N_{k+1}) , g_p^{-1}(N_{k+1})\big)} \\[7mm] 
\hspace{14mm} &=& \mu_p\big(f_p^{-1}(\ell^{k-1}_p)^{-1}(N_{k+1})) -  \mu_p(g_p^{-1}(N_{k+1})) 
\ < \ \sigma_p^{k-1}(N_{k+1}). 
\end{array}$
\vskip 5mm 
The continuity of the function $t = t(p)$ follows from the continuity of the functions $a_p(E_{B_k})$, $J^{\mu_p}\big(f_p^{-1}(\ell^{k-1}_p)^{-1}(B_k), g_p^{-1}(B_k)\big)$ in $p$ and   
$J^{\sigma^{k-1}_p}\big(\phi_t(B_k), B_k\big)$ in $(p,t)$ (cf. Lemma 3.2). 

These observations justify the definition of the map $\ell_k$ and it is readily seen to satisfy the required conditions. This completes the inductive step and 
we obtain the map $\ell^m$. 

The map $\ell^m$ satisfies the conditions on $\ell$ except $(1)'$.  
On the $n$-disk $A \cap L$ we compare the two maps 
$\sigma^m|_{A \cap L}, \tau|_{A \cap L} : P \to {\cal M}_g^\partial(A \cap L: m|_{A \cap L}\mbox{-reg})$, where $$\sigma^m = \ell^m_\ast((f_\ast\mu)|_A), \ \tau = g_\ast \mu : P \to 
{\cal M}_g^\partial(A: m|_{A}\mbox{-reg}).$$ 
Since \\
\hspace{18mm} 
$\begin{array}[t]{lcl}
\multicolumn{3}{l}{\ds \sigma^m(A \cap L) - \tau(A \cap L) 
\ = \ J^\mu(f^{-1}(\ell^m)^{-1}(A \cap L), g^{-1}(A \cap L))} \\[2mm] 
\hspace{15mm} &=& \ds J^\mu(f^{-1}(\ell^m)^{-1}(A), g^{-1}(A)) - \sum_{k=1}^m J^\mu(f^{-1}(\ell^m)^{-1}(B_k), g^{-1}(B_k)) \\[3mm] 
&=& \ds a(E_A) - \sum_{k=1}^m a(E_{B_k}) \ = \ 0, 
\end{array}$ \\[2mm] 
Theorem 2.1 yields a map 
\[ \xi : P \to {\cal H}_\partial(A \cap L; m|_{A \cap L} \mbox{-reg})_1 \cong {\cal H}_{\partial \cup (A - L)}(A; m|_A\mbox{-reg})_1 \] 
such that \ 
$(\xi_\ast \sigma^m)|_{A \cap L} = \tau|_{A \cap L}$ \ and \ 
$\xi_p = id_A$ if 
$\sigma_p|_{A \cap L} = \tau_p|_{A \cap L}$.
Finally the composition $\ell = \xi \,\ell^m$ satisfies all of the required conditions and this completes the proof. (We note that since the maps $\phi_t|_{N_k}$ ($t \in (-1,1)$) have locally common compact support, the map $h$ also has locally common compact support.) 
\end{proof} 


Let $L^0 = \emptyset$ and $f^0 \equiv id_{M_0}$, $g^0 \equiv id_{M_0}$. 

\begin{lemma}\label{l-sequence} 
There exists a sequence $(K_k, L_k, f^k, g^k)$ $(k=1, 2, \cdots)$ which satisfies the following conditions : 
\begin{itemize} 
\setlength{\itemsep}{1mm} 
\item[$(1_k)$] $K_k, L_k \in {\cal D}(M_0)$ and 
$L_{k-1}  \subset {\rm Int}_{M_0} K_k$, $K_k  \subset {\rm Int}_{M_0} L_k$
\item[$(2_k)$] 
\begin{itemize} 
\setlength{\itemsep}{1mm} 
\item[(i)\,] $f^k, g^k : P \to {\cal H}_{\partial}^c(M_0; m_0\mbox{-reg})_1^\ast$ are continuous maps
\item[(ii)] 
$f^k = \phi^k f^{k-1}$ and $g^k = \psi^k g^{k-1}$ for some continuous maps 
\[ \mbox{
$\phi^k : P \to {\cal H}_{\partial \cup L_{k-1}}^c(M_0; m_0\mbox{-reg})_1^\ast$ \ \ and \ \ 
$\psi^k : P \to {\cal H}_{\partial \cup K_k}^c(M_0; m_0\mbox{-reg})_1^\ast$
}\]
\end{itemize} 
\item[$(3_k)$] 
\begin{itemize} 
\setlength{\itemsep}{1mm} 
\item[(i)\,] 
$\ds {\rm diam}\,A \leq \frac{1}{2^k}$, \ 
$\ds {\rm diam}\,(g^{k-1}_p)^{-1}(A) \leq \frac{1}{2^k}$ \ $(A \in {\cal C}(cl_{M_0}(M_0-K_k)))$
\item[(ii)] 
$\ds {\rm diam}\,B \leq \frac{1}{2^k}$, \ 
$\ds {\rm diam}\,(f^k_p)^{-1}(B) \leq \frac{1}{2^k}$ \ $(B \in {\cal C}(cl_{M_0}(M_0 - L_k)))$
\end{itemize} 
\vskip 2mm 
\item[$(4_k)$] 
\begin{itemize} 
\setlength{\itemsep}{1mm} 
\item[(i)\ ] $f^k_\ast \mu  = g^{k-1}_\ast \mu$ on $K_k$ and 
$g^k_\ast \mu = f^k_\ast \mu$ on $L_k$ 
\item[(ii)\,] $J^\mu((f^k)^{-1}(A), (g^{k-1})^{-1}(A)) = a(E_A)$ \ $(A \in {\cal C}(cl_{M_0}(M_0-K_k))$
\item[(iii)] $J^\mu((f^k)^{-1}(B), (g^k)^{-1}(B)) = a(E_B)$ \ $(B \in {\cal C}(cl_{M_0}(M_0-L_k))$
\end{itemize} 
\item[$(5_k)$] 
\begin{itemize} 
\setlength{\itemsep}{1mm} 
\item[(i)\,] $\big\{(f^k_p)^{-1}\big\}_p$  is equi-continuous on $cl_{M_0}(M_0 - K_k)$ 
with respect to $d|_{M_0}$.
\item[(ii)] $\big\{(g^k_p)^{-1}\big\}_p$ is equi-continuous on $cl_{M_0}(M_0 - L_k)$ 
with respect to $d|_{M_0}$.
\end{itemize} 

\item[$(6_k)$] If $p \in P$ and $a_p = 0$, then $f^k_p = g^k_p = id_{M_0}$. 
\end{itemize} 
\end{lemma}

\begin{proof} Suppose we have constructed $(K_{k-1}, L_{k-1}, f^{k-1}, g^{k-1})$.  

Since $\big\{ (g^{k-1}_p)^{-1} \big\}_p$ is equicontinuous on $cl_{M_0}(M_0 - L_{k-1})$,  we can find $K_k \in {\cal D}(M_0)$ which satisfies $(1_k)$ and $(3_k)$. 
By applying Lemma \ref{l-basic} to the data ($L_{k-1}, K_k$, $f^{k-1}$, $g^{k-1}$, $\mu$, $a$), we obtain $\phi_k$ and 
$f^k$ which satisfies $(2_k)$, $(4_k)$ - $(6_k)$. 

Since $\big\{ (f^k_p)^{-1} \big\}_p$ is equicontinuous on $cl_{M_0}(M_0 - K_k)$, 
we can find $L_k \in {\cal D}(M_0)$ which satisfies $(1_k)$ and $(3_k)$. 
By applying Lemma \ref{l-basic} to the data ($K_k, L_k$, $g^{k-1}$, $f^k$, $\mu$, $-a$), we obtain $\psi^k$ and 
$g^k$ which satisfies $(2_k)$, $(4_k)$ - $(6_k)$. 
This completes the inductive step. 
\end{proof}

\begin{lemma}\label{l-limit}  Suppose $(K_k, L_k, f^k, g^k)$ $(k=1,2,\cdots)$ is the sequence in Lemma \ref{l-sequence} . 
\begin{itemize} 
\item[(1)] The sequence of maps $f^k : P \to {\cal H}_{\partial}(M_0; m_0\mbox{-reg})_1$  
$(k = 1,2, \cdots)$ converges $d|_{M_0}$-uniformly to a continuous map 
$f : P \to {\cal H}_{\partial}(M_0; m_0\mbox{-reg})_1$. 

\item[(2)] The sequence of maps $g^k : P \to {\cal H}_{\partial}(M_0; m_0\mbox{-reg})_1$  
$(k = 1,2, \cdots)$ converges $d|_{M_0}$-uniformly to a continuous map 
$g : P \to {\cal H}_{\partial}(M_0; m_0\mbox{-reg})_1$.  

\item[(3)] $f^{-1}|_{L_k} = (f^{k})^{-1}|_{L_k}$ and $g^{-1}|_{K_k} = (g^{k-1})^{-1}|_{K_k}$   
$(k = 1,2, \cdots)$ 

\item[(4)] $f_\ast \mu = g_\ast \mu$  
\item[(5)] If $p \in P$ and $a_p = 0$, then $f_p = g_p = id_{M_0}$.  
\end{itemize}
\end{lemma} 

\begin{proof} This follows from the same argument as in \cite[Proof of Lemma 4.8]{Be3}. 
\end{proof}

\begin{proof}[\bf Proof of Theorem 1.2$'$]  
We show that the continuous map 
$h = g^{-1} f : P \to {\cal H}_\partial(M_0, m_0\mbox{-reg})_1$, 
$h_p = g^{-1}_p f_p$, satisfies the required conditions.

(1) By Lemma \ref{l-limit}\,(4) we have $h_\ast \mu = \mu$ and from Lemma 2.2 it follows that 
\[ \mbox{$h_p \in {\cal H}_\partial(M_0, \mu_p\mbox{-reg})_1 \cap {\cal H}_\partial(M_0, \mu_p)
= {\cal H}_{\partial}(M_0, \mu_p)_1$. } \] 

(2) 
For each $F \in {\cal Q}(E_{M_0})$ there exists $k \geq 1$ and $A_1, \cdots, A_m \in {\cal C}(cl_{M_0}(M_0 - K_k))$ ($A_i \neq A_j$ ($i \neq j$)) such that $F = E_{A_1} \cup \cdots \cup E_{A_m}$ (disjoint). Thus, it suffices to show that $c^{\mu_p}_{h_p}(E_A) = a_p(E_A)$ for each $k \geq 1$ and $A \in {\cal C}(cl_{M_0}(M_0 - K_k))$. 

Since $f_p^{-1} \in {\cal H}_{E_{M_0}}(M_0)$,  we have $E_{f_p^{-1}(A)} = E_{A}$. 
Since 
$f_p^{-1}|_{K_k} = (f^{k}_p)^{-1}|_{K_k}$ and 
$g_p^{-1}|_{K_k} = (g^{k-1}_p)^{-1}|_{K_k}$ (Lemma 4.3\,(3)), we have 
\[ f_p^{-1}(A) = (f^k_p)^{-1}(A) \hspace{5mm} \text{and} \hspace{5mm} 
g_p^{-1}(A) = (g^{k-1}_p)^{-1}(A). \]  
Then from Lemma \ref{l-sequence} $(4_k)$ it follows that \\[2mm] 
\hspace{20mm} 
$\begin{array}[t]{lll} 
c^{\mu_p}_{h_p}(E_A) 
& = & c^{\mu_p}_{h_p}(E_{ f_p^{-1}(A)}) 
\ = \ J^{\mu_p}(f_p^{-1}(A), h_p f_p^{-1}(A)) 
\ = \ J^{\mu_p}(f_p^{-1}(A), g_p^{-1}(A)) \\[4mm] 
& = & J^{\mu_p}\big((f^k_p)^{-1}(A), (g^{k-1}_p)^{-1}(A)\big)
\ = \ a_p(E_A). 
\end{array}$
\vskip 3mm 
(3) From Lemma \ref{l-limit}\,(5) it follows that $h_p = id_{M_0}$. 
\end{proof}


\section{Proof of Theorem 1.2 in general case} 

In this final section we prove Theorem 1.2 in general case.
According to the usual strategy (cf. \cite{Br}),  
the mapping theorem in \cite{Be1, Be3} is used 
to reduce the noncompact $n$-manifold case to the $n$-cube with ends case (Theorem 1.2$'$). 
The correspondence between these cases under the proper map given by the mapping theorem has been discussed in Section 3.3. 

Throughout this sectioin $M^n$ is a noncompact connected  
$n$-manifold and $\omega \in {\cal M}_g^\partial(M)$. 

\begin{lemma}\label{l-mapping} $($\cite[Proposition 4.2, Proof of Theorem 4.1 (p\,252)]{Be3}$)$  
There exists 
a compact 0-dimensional subset $E \subset \partial I^n$ $(E \subset I_1$ if $n \geq 2)$ and 
a continuous proper surjection $\pi : I^n - E \to M$ 
which satisfies the following conditions: 
\begin{itemize}
\item[(i)\ ] $U \equiv \pi({\rm Int}\,I^n)$ is a dense open subset of ${\rm Int}\,M$ and 
$\pi|_{{\rm Int}\,I^n} : {\rm Int}\,I^n \to U$ is a homeomorphism.
\item[(ii)\,] $F \equiv \pi(\partial I^n - E) = M - U$ and $\omega(F) = 0$. 
\item[(iii)]  The induced map $\overline{\pi} : E \to E_M$ is a homeomorphism. 
\item[(iv)] The induced measure $\pi^\ast \omega$ is $m|_{I^n - E}$-regular. 
\end{itemize} 
\end{lemma} 

Let $M_0 = I^n - E$ and $m_0 = m|_{M_0}$. We have 
$\omega_0 \equiv \pi^\ast \omega \in {\cal M}^\partial_g(M_0, m_0\mbox{-reg})$. 


\begin{proof}[\bf Proof of Theorem 1.2] 
By the considerations in Section 3.3 the map $\pi$ in Lemma \ref{l-mapping} induces 
the reciprocal homeomorphisms in the left side and the commutative diagram of three squares: \\[2mm] 
\hspace{2mm} 
$\begin{array}[t]{clcccccccl}
& & & & & c^{\pi^\ast \mu_p} & & & & \\[-1mm] 
{\cal M}_g^\partial(M_0, m_0\mbox{-reg}) & \hspace{5mm} & {\cal H}_\partial(M_0, m_0\mbox{-reg})_1 & \supset & {\cal H}_\partial(M_0, \pi^\ast\mu_p)_1  
& \lra & {\cal S}(M_0, \pi^\ast \mu_p) & \subset & {\cal S}(E_{M_0}) & \\[2mm] 
\pi_\ast \ \Big\downarrow \Big\uparrow \ \pi^\ast 
& 
& \Big\downarrow & \hspace{-46mm} \pi_\ast & \Big\downarrow & \hspace{-47.5mm} \pi_\ast & \Big\downarrow & \hspace{-29.5mm} \cong \hspace{3.5mm} \overline{\pi}_\ast & \Big\downarrow & \hspace{-15.5mm} \cong \hspace{4mm} \overline{\pi}_\ast \\[3mm]  
{\cal M}_g^\partial(M, \omega\mbox{-reg}) & & {\cal H}_F(M, \omega\mbox{-reg})_1 & \supset & {\cal H}_F(M, \mu_p)_1 & \lra & {\cal S}(M, \mu_p).  & \subset & {\cal S}(E_{M}) & \\[-0.5mm]
& & & & & c^{\mu_p} & & & & \\[-2mm]
\,\pi^\ast = (\pi_\ast)^{-1} & & & & & & & & &
\end{array}$ 
\vskip 4mm 
The maps $\mu$ and $a$ admit the lifts to $M_0$ : 
\[ \mbox{$\pi^\ast \mu : P \to {\cal M}_g^\partial(M_0, m_0\mbox{-reg})$ \ \ and \ \ 
$\widetilde{a} = (\overline{\pi}_\ast)^{-1} a : P \to {\cal S}(E_{M_0})$. } \] 
Since $a_p \in S(M, \mu_p)$, the 3rd square in the above diagram implies $\widetilde{a}_p \in 
{\cal S}(M_0, \pi^\ast \mu_p)$. 
Theorem 1.2$'$ provides with a continuous map 
$\widetilde{h} : P \to {\cal H}_\partial(M_0, \pi^\ast \omega\mbox{-reg})_1$ 
such that for each $p \in P$ \\[1.5mm]  
\hspace{5mm} 
(1)$'$ $\widetilde{h}_p \in {\cal H}_{\partial}(M_0, \pi^\ast \mu_p)_1$, 
\hspace{5mm} 
(2)$'$ $c^{\pi^\ast \mu_p}_{\widetilde{h}_p} = \widetilde{a}_p$, 
\hspace{5mm} 
(3)$'$ if $\widetilde{a}_p = 0$, then $\widetilde{h}_p = id_{M_0}$. 
\vskip 1.5mm 
We show that the map \hspace{2mm} 
\[ \mbox{$h = \pi_\ast \widetilde{h} : P \to {\cal H}_F(M, \omega\mbox{-reg})_1 \subset {\cal H}_\partial(M, \omega\mbox{-reg})_1$} \] 
satisfies the required conditions.  
\begin{itemize}
\item[(1)] The condition (1)$'$ and the 1st square imply that $h_p \in {\cal H}_F(M, \mu_p)_1 \subset {\cal H}_\partial(M, \mu_p)_1$. 

\item[(2)] From (2)$'$ and the 2nd square it follows that 
\[ c^{\mu_p}_{h_p} 
= c^{\mu_p} \pi_\ast \big(\widetilde{h}_p \big) 
= \overline{\pi}_\ast c^{\pi^\ast \mu_p}\big (\widetilde{h}_p \big)
= \overline{\pi}_\ast \big(c^{\pi^\ast \mu_p}_{\widetilde{h}_p}\big)
= \overline{\pi}_\ast \big(\widetilde{a}_p \big)
= a_p.
\]
\vskip 2mm 

\item[(3)] If $a_p = 0$, then $\widetilde{a}_p = 0$ and $\widetilde{h}_p = id_{M_0}$. 
This implies that $h_p = id_M$.  
\end{itemize} 
This completes the proof. 
\end{proof}

\begin{proof}[\bf Proof of Theorem 1.1] 
The required section is obtained by applying Theorem 1.2 to the data: 
$P = {\cal S}(M, \omega)$, $\mu \equiv \omega$ and $a$ is the inclusion ${\cal S}(M, \omega) \subset {\cal S}(E_M)$. 
\end{proof} 

Suppose ${\cal G}$ is any subgroup of ${\cal H}_{E_M}(M, \omega)$ with 
${\cal H}_\partial(M, \omega)_1 \subset {\cal G}$. 
Consider the restriction \break $c^\omega|_{\cal G} : {\cal G} \to {\cal S}(M, \omega)$. 

\begin{corollary} {\rm (1)} $({\cal G}, {\rm Ker}\,c^\omega|_{\cal G}) \cong ({\rm Ker}\,c^\omega|_{\cal G}) \times ({\cal S}(M, \omega), 0)$. 
\begin{itemize}
\item[(2)] ${\rm Ker}\,c^\omega|_{\cal G}$ is a strong deformation retract of ${\cal G}$.
\end{itemize} 
\end{corollary} 

\begin{proof} 
(1) The required homeomorphism is defined by 
\[ \mbox{$\phi : {\cal G} \to ({\rm Ker}\,c^\omega|_{\cal G}) \times {\cal S}(M, \omega)$, \ \ 
$\phi(h) = ((s(c^\omega_h))^{-1}h, c^\omega_h)$.} \] 
The inverse is given by $\phi^{-1}(f, a) = s(a)f$. 

(2) Since the topological vector space ${\cal S}(M, \omega)$ admits a strong deformation retraction onto $\{ 0 \}$, 
the conclusion follows from (1). 
\end{proof}


\end{document}